\documentclass[letterpaper,twoside,draft,11pt]{article}
\usepackage{amsmath}               
\usepackage{amsthm}                

\usepackage{latexsym}

\usepackage{amssymb}
\usepackage{xspace}
\usepackage{amscd}

\date{}

\newcommand{\Gal}{{\rm Gal}}

\newcommand{\Imagen}{\mbox{\rm Im}}

\newcommand{\End}{\mbox{\rm End}}

\newcommand{\inv}{^{-1}}

\newcommand{\N}{{\mathbb N}}
\newcommand{\Z}{{\mathbb Z}}
\newcommand{\Q}{{\mathbb Q}}

\newcommand{\C}{{\mathbb C}}

\newcommand{\tr}{{\rm tr}}
\newcommand{\diag}{{\rm diag}}
\newcommand{\matriz}[1]{\begin{array} #1 \end{array}}
\newcommand{\pmatriz}[1]{\left(\begin{array} #1 \end{array}\right)}
\newcommand{\GEN}[1]{\langle #1 \rangle}

\title{Bicyclic units, Bass cyclic units and free groups\thanks{
Research supported by Capes of Brazil, CNPq grant 303.756/82-5 and Fapesp, Projeto Tematico 00/07.291-0, D.G.I. of
Spain and Fundaci\'{o}n S\'{e}neca of Murcia.}}

\author{\begin{tabular}{cc}
Jairo Z. Gon\c{c}alves & \'{A}ngel del R\'{i}o \\
Departamento de Matem\'{a}tica & Departamento de Matem\'{a}ticas \\
Universidade de S\~ao Paulo & Universidad de Murcia \\
S\~ao Paulo, 05389-970, Brazil & Campus de Espinardo, Murcia 30100, Spain\\
jzg@ime.usp.br & adelrio@um.es
\end{tabular}}

\newtheorem{theorem}{Theorem}[section]
\newtheorem{lemma}[theorem]{Lemma}
\newtheorem{proposition}[theorem]{Proposition}

\newtheorem{corollary}[theorem]{Corollary}

\theoremstyle{remark}\newtheorem{remark}[theorem]{Remark}

\theoremstyle{remark}\newtheorem{example}[theorem]{Example}

\theoremstyle{remark}

\theoremstyle{remark}\newtheorem{question}[theorem]{Question}

\theoremstyle{remark}\newtheorem{problem}[theorem]{Problem}

\begin{document}

\maketitle

\begin{abstract}
Let $G$ be a finite group and $\Z G$ its integral group ring. We show that if $\alpha$ is a non-trivial bicyclic unit
of $\Z G$, then there are bicyclic units $\beta$ and $\gamma$ of different types, such that $\GEN{\alpha,\beta}$ and
$\GEN{\alpha,\gamma}$ are non-abelian free groups. In case that $G$ is non-abelian of order coprime with $6$, then we
prove the existence of a bicyclic unit $u$ and a Bass cyclic unit $v$ in $\Z G$, such that for every positive integer
$m$ big enough, $\GEN{u^m,v}$ is a free non-abelian group.
\end{abstract}

\section{Introduction}

A {\em free pair} is by definition a pair formed by two generators of a non-abelian free group. Let $G$ be a finite
group. The existence of free pairs in the group of units $U(\Z G)$ of the integral group ring $\Z G$, was firstly
proved by Hartley and Pickel \cite{HP}, provided that $G$ is neither abelian nor a Hamiltonian $2$-group (equivalently
$U(\Z G)$ is neither abelian nor finite). Their proof is not constructive and this raised the question of exhibiting a
concrete free pair. This goal was achieved for non-Hamiltonian groups by Marciniak and Sehgal \cite{MS} using bicyclic
units, and for Hamiltonian groups (non 2-group) by Ferraz \cite{F} using Bass cyclic units. These results, together
with a classical theorem of Bass \cite{Bass}, that states that if $G$ is abelian then the Bass cyclic units generates a
subgroups of finite index in $U(\Z G)$, and the more recent ones of Ritter and Sehgal \cite{RS} and Jespers and Leal
\cite{JL}, which prove that the group generated by the bicyclic and the Bass cyclic units generates a big portion of
$U(\Z G)$, show that these two types of units have an important role in the structure of $U(\Z G)$. As a consequence,
several authors have payed attention to the problem of describing the structure of the group generated either by
bicyclic units, or Bass cyclic units and more specifically, to the problem of deciding when two bicyclic units or Bass
cyclic units form a free pair \cite{DJR,GPJA,JRR,S}.

The {\em bicyclic units} of $\Z G$ are the elements of one of the following forms
    $$\beta_{x,h}=1+(1-h)x(1+h+h^2+\ldots+h^{d-1}), \quad
    \gamma_{x,h}=1+(1+h+h^2+\ldots+h^{d-1})x(1-h),$$
where $x,h\in G$ and $h$ has order $d$. Two bicyclic units are of {\em same type} if both are of type $\beta$ or both
are of type $\gamma$. Otherwise they are of {\em different types}.

The Bass cyclic units of $\Z G$ are the elements of the form
    $$u_{k,m}(x) = (1+x+x^2+\ldots+x^{k-1})^m +\frac{1-k^m}{d}(1+x+x^2+\ldots+x^{d-1}),$$
where $x$ is an element of order $d$ and $k$ and $m$ are positive integers such that $\gcd(k,d)=1$ and $\phi(d)|m$.

One says that a complex number $z$ is {\em free} or that $z$ is a {\em free point}, if the matrices
    $$
    \pmatriz{{cc} 1 & 1 \\ 0 & 1} \quad \pmatriz{{cc} 1 & 0 \\ z & 1}
    $$
form a free pair\footnote{Some authors define free points by considering $\pmatriz{{cc} 1 & 2 \\ 0 & 1}$ as the first
matrix. The two definitions are somehow equivalent because $z$ is free under this definition if and only if $2z$ is
free under our definition.}. Otherwise $z$ is said to be {\em non-free}. The most classical result on the subject is
Sanov Theorem \cite{Sv} which states that every complex number of modulus at least $4$ is free. On the other hand every
integer of modulus less than $4$ is non-free. For an up-to-date list of results on the matter see \cite{Ba}.

If $R$ is a ring of characteristic $0$ then its elements can be considered as belonging to a complex algebra. Indeed,
$R$ embeds in $R\otimes_{\Z} \C$, since $_Z\C$ is torsion-free and so it is flat. Thus it makes sense to talk on the
transcendency of these elements over any subfield of $\C$. Moreover, left multiplication by elements of $R$ can be
considered as endomorphisms of the underlying vector space structure of the algebra. So we can refer to the eigenvalues
of such endomorphism. Note that the eigenvalues do not depend on the algebra where $R$ is included because the minimal
polynomial is independent of the algebra. The first result of this paper is the following freeness criteria in terms of
free points.

\begin{theorem}\label{Criteria}
Let $R$ be a ring of characteristic $0$ and $a$ and $b$ elements of $R$ such that $a^2=b^2=0$. Then $(1+a,1+b)$ is a
free pair if and only, one of the following conditions hold:
\begin{enumerate}
\item $ab$ is transcendental (over the rationals).
\item $ab$ is algebraic (over the rationals) and one of the eigenvalues of $ab$ is free.
\end{enumerate}
\end{theorem}

Notice that every bicyclic unit is of the form $1+a$ for some $a\in \Z G$ such that $a^2=0$ and so
Theorem~\ref{Criteria} can be applied to pairs of bicyclic units. As an application of Theorem~\ref{Criteria} we show a
method to construct many free pairs of bicyclic units and in particular obtain.

\begin{theorem}\label{FreeComp}
If $\Z G$ has a non-trivial unit $\alpha$ then $\Z G$ has bicyclic units $\beta$ and $\gamma$ of different type, such
that $(\alpha,\beta)$ and $(\alpha,\gamma)$ are free pairs.
\end{theorem}

In fact the bicyclic units $\beta$ and $\gamma$ of Theorem~\ref{FreeComp} can be explicitly constructed (see
Corollary~\ref{FreeCompExp}). The existence of a free pair formed by bicyclic units of different types was already
proven in \cite{MS}. Another consequence of Theorem~\ref{Criteria} is a result of Salwa, which states that if $a$ and
$b$ satisfy the conditions of Theorem~\ref{Criteria} and $ab$ is non-nilpotent then $(1+a,(1+b)^m=1+mb)$ is a free pair
for some positive integer $m$. In Section~\ref{ExBic} we discuss the minimal $m$ for which $(1+a,1+mb)$ is a free pair.

Then we consider groups generated by a bicyclic unit and a Bass cyclic unit. Gon\c{c}alves and Passman \cite{GPJA} proved
recently that if $G$ is a finite non-abelian group of order coprime with 6, then $\Z G$ has a free pair formed by Bass
cyclic units. With the same assumptions we also prove:

\begin{theorem}\label{BB}
If $G$ is a finite non-abelian group of order coprime with 6, then $\Z G$ has a bicyclic unit $\beta$ and a Bass cyclic
unit $u$ such that $(u,\beta^t)$ is a free pair for any sufficiently large positive integer $t$.
\end{theorem}

The hypothesis of the order of $G$ being coprime with $6$, in Theorem~\ref{BB}, is partially justified because for the
existence of a non-trivial Bass cyclic unit one needs the exponent of the group $G$ not being a divisor of $4$ or $6$.
In the last section of the paper we prove that $\Z S_n$ (respectively, $\Z A_n$) has a free pair formed by a Bass
cyclic unit and a bicyclic unit if and only if $n\ge 5$.

\section{Proof of Theorem~\ref{Criteria}}\label{SecCrit}

For the proof of Theorem~\ref{Criteria} we first prove some lemmas. The first one was already proved in \cite{GPUnp}
for $K=\C$ and the same proof works in general. We include a proof for completeness.

\begin{lemma}\label{TransLem}
Let $K$ be a field of characteristic $0$ and $a$ and $b$ be elements of a $K$-algebra such that $a^2=b^2=0$. If $ab$ is
transcendental over $K$ then $K[a,b]$ is naturally isomorphic to the relatively free algebra $K[x,y|x^2=y^2=0]$.
\end{lemma}

\begin{proof}
Let $f_1,f_2,f_3$ and $f_4$ be polynomials in one variable with coefficients in $K$ such that
    $$f_1(ab)+f_2(ab)a+bf_3(ab)+bf_4(ab)a = 0.$$
Multiplying on both sides by $ab$ one has $abf_1(ab)ab=0$ and hence $f_1=0$, by the transcendency of $ab$ over $K$.
Then multiplying by $ab$ on the left and by $b$ on the right, one has $abf_2(ab)ab=0$ and so $f_2=0$. By symmetry one
has $f_3=0$. Finally multiplying by $a$ on the left and by $b$ on the right one deduces $abf_4(ab)ab=0$ and this yields
$f_4=0$.

This proves that the elements of the form $(ab)^i, (ab)^ia, b(ab)^i$ and $b(ab)^i a$ are linearly independent over $K$.
Therefore the algebra homomorphism $f:K[x,y|x^2=y^2=0] \rightarrow K[a,b]$ given by $f(x)=a$ and $f(y)=b$ is an
isomorphism.
\end{proof}

\begin{lemma}\label{Producto}
Let $R=R_1\times \cdots \times R_n$ be a direct product of rings, and $u=(u_1,\ldots,u_n)$ and $v=(v_1,\ldots,v_n)$ be
units of $R$, with $u_i,v_i\in R_i$ for every $i$. Then $(u,v)$ is a free pair if and only if $(u_i,v_i)$ is a free
pair for some $i$.
\end{lemma}

\begin{proof}
Without loss of generality we can assume that $n=2$. The sufficient condition is clear. Assume that neither $(u_1,v_1)$
nor $(u_2,v_2)$ is a free pair. Then there are non-trivial words $w_1$ and $w_2$ in the free group on two symbols such
that $a_i=w_i(u_i,v_i)=1$ for $i=1,2$. If the commutator $w=w_1w_2w_1\inv w_2\inv \ne 1$ then
    $$\matriz{{rcl}
    w(u,v)&=& (a_1 w_2(u_1,v_1) a_1^{-1} w_2^{-1}(u_1,v_1),w_1(u_2,v_2) a_2 w_1^{-1}(u_2,v_2) a_2^{-1}) \\
    &=& ((w_2(u_1,v_1)w_2^{-1}(u_1,v_1),w_1(u_2,v_2)w_1^{-1}(u_2,v_2)) = (1,1)=1.}$$
Otherwise $w_1$ and $w_2$ belongs to a cyclic group and therefore there is  $1\ne w \in \GEN{w_1}\cap \GEN{w_2}$. Then
$w(u,v)=1$. So in both cases $(u,v)$ is not a free pair.
\end{proof}

%

\begin{lemma}\label{Dos}
If $R = M_m(\C)=\C[a,b]$ with $a^2=b^2=0$, then $m\le 2$. Furthermore $m=1$ if and only if $ab$ is nilpotent.
\end{lemma}

\begin{proof}
For a real number $\alpha$, let $[\alpha]$ denote the greatest integer non greater than $\alpha$. Consider the elements
of $R$ as endomorphisms of $\C^n$. Then $\dim_{\C} \Imagen(ab) \le \Imagen(a) \le \left[\frac{m}{2}\right]$, because
$\Imagen(a) \subseteq \ker(a)$. Then, by the Cayley-Hamilton Theorem, the degree of the minimal polynomial of $ab$ over
$\C$ is $\le \left[\frac{m}{2}\right]+1$. On the other hand $R=\C[ab]+b\C[ab]+\C[ab]a+b\C[ab]a$. Hence $m^2 = \dim_{\C}
\End_{\C}(\C^m) \le 4\left(\left[\frac{m}{2}\right]+1\right)$ and so $m\le 2$. This shows that $m\le 2$.

If $m=1$ then clearly $ab$ is nilpotent. By means of contradiction assume that $m=2$ and $ab$ is nilpotent. After
conjugating by some invertible matrix one may assume that $a=\pmatriz{{cc}0&1\\0&0}$. Then the second row of $ab$ is
$0$ and, since $ab$ is nilpotent, also the $(1,1)$-entry is $0$. Then $b=\pmatriz{{cc} 0&\mu\\0&0}$ for some $\mu$.
Therefore $ab=ba$, contradicting the assumption $M_2(\C)=\C[a,b]$.
\end{proof}

{\em Now we prove Theorem~\ref{Criteria}}.

Let $R$ be a ring of characteristic $0$ and $a,b\in R$ such that $a^2=b^2=0$. By replacing $R$ by $\C\otimes_{\Z} R$ if
needed one may assume that $R$ is a $\C$-algebra. If $ab=ba$ then $\GEN{1+a,1+b}$ is not free and $0$ is the only
eigenvalue of $ab$, hence the Theorem holds. In the remainder of the proof we assume that $ab\ne ba$.

Assume first that $R=M_2(\C)=\C[a,b]$ and identify $R$ with the ring of endomorphism of a $2$-dimensional vector space
over $\C$. Since $\ker(ab)\subseteq \ker(b) \ne 0$, one of the eigenvalues of $ab$ is $0$. Let $\lambda$ be the other
eigenvalue. By Lemma~\ref{Dos}, $ab$ is not nilpotent and so $\lambda\ne 0$. Hence $ab$ is diagonalizable and after
some suitable conjugation, one may assume that $ab=\pmatriz{{cc} \lambda&0\\0&0}$. Then $\ker(b) = \ker(ab)$ and
$\Imagen(a)= \Imagen(ab)$ and hence
    $$a=\pmatriz{{cc} 0 & \mu \\ 0 & 0} \quad \mbox{ and } \quad b=\pmatriz{{cc} 0 & 0 \\ \frac{\lambda}{\mu} & 0}$$
for some $0\ne \mu\in \C$. Conjugating by the matrix $\pmatriz{{cc} \frac{1}{\mu} & 0 \\
0 & 1}$ one may assume that
    $$a=\pmatriz{{cc} 0 & 1 \\ 0 & 0} \quad \mbox{and} \quad b=\pmatriz{{cc} 0 & 0 \\ \lambda & 0}$$
and so $(1+a,1+b)$ is a free pair if and only if $\lambda$ is free. Note that $ab$ is transcendental over $\Q$ if and
only if so is $\lambda$ and in this case $\lambda$ is free.

Now we consider the general case.

If $ab$ is transcendental (over $\Q$) then, by Lemma~\ref{TransLem}, there is an algebra homomorphism
$\Q[a,b]\rightarrow M_2(\Q)$ mapping $a$ to $\alpha=\pmatriz{{cc} 0 & 1 \\ 0 & 0}$ and $b$ to $\beta=\pmatriz{{cc} 0 &
0 \\ 4 & 0}$. Since $4$ is a free-point, $(1+\alpha,1+\beta)$ is a free pair and thus $(1+a,1+b)$ is also a free pair.

Assume now that $ab$ is algebraic (over $\Q$). We may assume without loss of generality that $R=\C[a,b]$. Then
$\dim_{\C} \C[ab]<\infty$ and $R=\C[ab]+b\C[ab]+a\C[ba]+b\C[ab]a$. Thus $\dim_{\C} R\le 4 \dim_{\C} \C[ab]<\infty$.
Therefore, if $J$ denotes the Jacobson radical of $R$ then $J$ is nilpotent and $R/J$ is semisimple. The first implies
that $1+J$ is a nilpotent normal subgroup of the group of units of $R$. Hence $(1+J)\cap \GEN{1+a,1+b}$ is nilpotent
and so $\GEN{1+a,1+b}$ is free if and only if $\GEN{1+a+J,1+b+J}$ is free. Also from the nilpotency of $J$ it follows
that the minimal polynomial of $ab$ divides a power of the minimal polynomial of $ab+J$. We conclude that one may
assume without loss of generality that $J=0$ and so $R$ is semisimple. Then $R= \oplus_{i=1}^k A_i$ with $A_i\simeq
M_{n_i}(\C)$ and $n_i\le 2$ for every $i$, by Lemma~\ref{Dos}. If $x\in R$ and $1\le i \le k$ then $x_i$ denotes the
projection of $x$ in $A_i$.

Assume first that one of the eigenvalues $\lambda$ of $ab$ is free. Then $\lambda$ is an eigenvalue of $a_ib_i$ for
some $1\le i \le k$. Moreover, $n_i=2$ because if $n_i=1$ then $a_i=b_i=0$. By the first part of the proof,
$(1+a_i,1+b_i)$ is a free pair of $A_i$. Thus $\GEN{1+a,1+b}$ is a free pair.

Conversely, assume that $(1+a,1+b)$ is a free pair of $R$. By Lemma~\ref{Producto}, there is $1\le i \le k$ such that
$(1+a_i,1+b_i)$ is a free pair. Clearly $n_i=2$ and, by the first part of the proof, one of the eigenvalues of $a_ib_i$
is free. Then one of the eigenvalues of $ab$ is free and this finishes the proof of Theorem~\ref{Criteria}.

\section{Applications of Theorem~\ref{Criteria}}\label{SecAppl}

A first application of Theorem~\ref{Criteria} is the following corollary. The second statement appeared in \cite{S} and
a weak version of the third one appeared in \cite{JRR}.

\begin{corollary}\label{Nil}
Let $R$ be a ring of characteristic $0$. Let $\alpha=1+a$ and $\beta=1+b$ be units of $R$ with $a^2=b^2=0$.
\begin{enumerate}
\item If $m$ and $n$ are positive integers then $(\alpha^m,\beta^n)$ is a free pair if and only if
    $(\alpha,\beta^{mn})$ is a free pair.
\item If $ab$ is not nilpotent then $(\alpha,\beta^m)$ is a free pair for some positive integer $m$.
\item $ab$ is nilpotent if and only if $\GEN{1+a,1+b}$ is nilpotent. Moreover, if $ab$ is algebraic over $\C$
then, $\GEN{1+a,1+b}$ is nilpotent if and only if $0$ is the only eigenvalue of $ab$.
\end{enumerate}
\end{corollary}

\begin{proof}
(1) Is an obvious consequence of Theorem~\ref{Criteria} because $\alpha^m=1+ma$ and $\beta^n=1+nb$.

(2) and (3). Firstly assume that $ab$ is nilpotent and let $S$ be the multiplicative semigroup generated by $a$ and
$b$. Then there is a positive number $n$ such that $s^n=0$ for each $s\in S$. This implies that $1+S$ is a nilpotent
subgroup of the group of units of $R$. Hence $\GEN{1+a,1+b}$ is nilpotent.

Secondly assume that $ab$ is transcendental over $\Q$. Then $(\alpha,\beta)$ is a free pair, by Theorem~\ref{Criteria},
and so (2) holds for $m=1$ and $\GEN{1+a,1+b}$ is not nilpotent.

Thirdly assume that $ab$ has a non-zero eigenvalue $\lambda$. Let $m$ be a positive integer such that $|m\lambda|\ge
4$. Then $m\lambda$ is an eigenvalue of $mab$ and $m\lambda$ is free by Sanov Theorem. Therefore
$(\alpha=1+a,\beta^m=1+mb)$ is a free pair, by Theorem~\ref{Criteria}, and hence $\GEN{\alpha,\beta}$ is not nilpotent.

Finally, if $ab$ is algebraic over $\Q$ and $0$ is the only eigenvalue of $ab$ then $ab$ is nilpotent.
\end{proof}

\begin{corollary}\label{Rep}
Let $a$ and $b$ be elements of a finite dimensional $\C$-algebra $A$, such that $a^2=b^2=0$. Then $(1+a,1+b)$ is a free
pair (resp. $\GEN{1+a,1+b}$ is nilpotent) if and only if there is an irreducible representation $\rho$ of $A$ such that
one of the eigenvalues of $\rho(ab)$ is free (resp. $0$ is the only eigenvalue of $\rho(ab)$ for each irreducible
representation $\rho$ of $A$).
\end{corollary}

\begin{proof}
It follows from Theorem~\ref{Criteria} because the set of eigenvalues of an element $x$ of $A$ is the union of the sets
of eigenvalues of $\rho(x)$, for $\rho$ running on the irreducible representations of $A$.
\end{proof}

\begin{corollary}\label{Char}
Let $G$ be a finite group, with the property that all the (complex) irreducible characters of $G$ have degree $\le 3$.
Let $\C G$ be the complex group algebra of $G$, and let $a,b\in \C G$ be such that $a^2=b^2=0$. Then $(1+a,1+b)$ is a
free pair (resp. $\GEN{1+a,1+b}$ is nilpotent) if and only if $\chi(ab)$ is free for some irreducible character $\chi$
of $G$ (resp. $\chi(ab)=0$ for every irreducible character $\chi$ of $G$).
\end{corollary}

\begin{proof}
Let $\rho$ be an irreducible representation of $G$ of degree $n$, $[n/2]$ the greater integer less of equal than $n/2$,
and $\chi$ the character afforded by $\rho$. Then $\dim_{\C} \Imagen(\rho(ab)) \le \dim_{\C} \Imagen(\rho(a)) \le [n/2]
\le 1$, because $\Imagen(\rho(a)) \subseteq \ker(\rho(a))$ and $n\le 3$. This shows that the eigenvalues of $\rho(ab)$
are $0$ and $\chi(ab)$. Now the result follows from Corollary~\ref{Rep}.
\end{proof}

A {\em trace map} on a complex algebra $A$ is a $\C$-linear form $T:A\rightarrow \C$ satisfying the following
conditions for $x,y\in A$: $T(xy)=T(yx)$; if $x$ is nilpotent then $T(x)=0$; and if $0\ne x=x^2$ then $T(x)$ is a
positive real number. For example, the classical trace $\tr:M_n(\C)\rightarrow \C$ is a trace map on $M_n(\C)$. An easy
argument shows that if $T$ is a trace map on $M_n(\C)$ then $T=\alpha \tr$ for some positive real number $\alpha$
(namely $\alpha=T(1)/n$).

\begin{corollary}\cite{S}\label{Salwa}
Let $R$ be a ring of characteristic $0$ and $a,b\in R$ with $a^2=b^2=0$. If $\C\otimes_{\Z} R$ has a trace map $T$ such
that $|T(ab)|\ge 2 T(1)$ then $(1+a,1+b)$ is a free pair.
\end{corollary}

\begin{proof}
By the definition of trace map, $ab$ is not nilpotent. If $ab$ is transcendental over $\Q$, then the result follows at
once from Theorem~\ref{Criteria}.

Otherwise, one may assume that $R=\C[a,b]$. By Lemma~\ref{Dos}, $R/J\simeq \prod_{i=1}^k M_{n_i}(\C)$ with $n_i\le 2$.
Suppose that $n_i=1$ if $i>l$ and $n_i=2$, otherwise. Since $R$ is artinian, $T(J)=0$, because $J$ is nilpotent, and
idempotents lift modulo $J$. Thus $T$ induces a trace map $\overline{T}$ on $R/J$. For each $i$, the restriction of
$\overline{T}$ to $M_{n_i}(\C)$ is a trace map $T_i=\alpha_i \tr$ on $M_{n_i}(\C)$. One of the eigenvalues of $a_ib_i$
is $0$. If $n_i=2$, then let $\lambda_i$ be the other eigenvalue of $a_ib_i$. Then
    $$4\sum_{i=1}^l \alpha_i \le 2\sum_{i=1}^k n_i \alpha_i = 2 T(1) \le T(a b)=\sum_{i=1}^k \alpha_i \tr(a_i b_i) =
        \sum_{i=1}^l \alpha_i \lambda_i$$
and therefore $|\lambda_i|\ge 4$ for some $i$. Since $\lambda_i$ is an eigenvalues of $ab$, $(1+a,1+b)$ is a free pair
by Theorem~\ref{Criteria}.
\end{proof}

\section{Bicyclic units}

In this section $G$ stands for an arbitrary group. The notation $H\le G$ means that $H$ is a subgroup of $G$. If $A$ is
a subset of $G$ then $\GEN{A}$ denotes the subgroup generated by $A$. If moreover $A$ is finite then we write
$\overline{A} = \sum_{a\in A} a\in \Z G$. If $g\in G$ has finite order then we abbreviate
$\overline{g}=\overline{\GEN{g}}$.

Let $H$ be a finite subgroup of $G$, $h\in H$ and $x\in G$. Then $(1-h)x\overline{H}$ and $\overline{H}x(1-h)$ are
elements of $\Z G$ of square zero and therefore
    $$\beta_{x,h,H} = 1+(1-h)x\overline{H} \quad \mbox{and} \quad \gamma_{x,h,H} = 1+\overline{H}x(1-h)$$
are units of $\Z G$. So the bicyclic units of $\Z G$ are the elements of one of the following forms
    $$\beta_{x,h}=\beta_{x,h,\GEN{h}}=1+(1-h)x\overline{h} \quad \mbox{and} \quad
    \gamma_{x,h}=\gamma_{x,h,\GEN{h}}= 1+\overline{h}x(1-h),$$
where $x,h\in G$ and $h$ has finite order.

If $a=\sum_{g\in G} a_g g \in \C G$ then the {\em trace} of $a$ is by definition $T(a)=a_1$. Notice that $T$ is a trace
map on $\C G$ as defined in Section~\ref{SecAppl}. This is clear if $G$ is finite because then the trace of the regular
representation of $\C G$ is $|G|T$. For a proof for infinite groups see \cite[Lemma 1.7 in page 37 and Lemma 3.3 in
page 47]{P}.

\begin{lemma}\label{trace} Let $x,y\in G$ and $H,K\le G$.
Then $\overline{H}x\overline{K} = |H\cap xKx\inv| \cdot \overline{HxK}$ and
    $$T(y\overline{H}x\overline{K}) = \left\{\matriz{{ll} |H\cap xKx\inv|, & \mbox{if } y\inv \in HxK \\ 0, &
    \mbox{otherwise.}}\right.$$
\end{lemma}

\begin{proof}
Consider the map $f:H\times K \rightarrow HxK$, given by $f(h,k)=hxk$. Clearly $\overline{H}x\overline{K} = \sum_{y\in
HxK} |f\inv(y)| y$. It is easy to see that the map $f\inv(x)\rightarrow f\inv(hxk)$, given by $(h_1,k_1)\mapsto
(hh_1,k_1k)$ is bijective. Using the equality $|H||K|=|HxK|\cap |H\cap xKx\inv|$ one deduces that $|f\inv(y)|=|H\cap
xKx\inv|$ for every $y\in HxK$ and hence $\overline{H}x\overline{K}=|H\cap xKx\inv|\overline{HxK}$.
%
%
Then $T(y\overline{H}x\overline{K}) = 0$ if $y\inv \not\in HxK$ and otherwise $T(y\overline{H}x\overline{K}) = |H\cap
xKx\inv|$, as asserted.
\end{proof}

\begin{proposition}\label{ManyFP}
Let $H\le K$ be finite subgroups of a group $G$ and let $x\in G$, $h\in H$ and $k\in K$ be such that $x\inv h x\not\in
K$.
\begin{enumerate}
\item If $x\inv k x\not\in K$ then $(\beta_{x,h,H},\gamma_{x\inv,k,K})$ is a free pair.
\item If $xk x\inv \not\in K$ then $(\beta_{x,h,H},\beta_{x\inv,xkx\inv,xKx\inv})$ is a free pair.
\end{enumerate}
\end{proposition}

\begin{proof}
Let $a=\gamma_{x\inv,k,K}-1$ and $b=\beta_{x,h,H}-1$. Then
    $$ab = \overline{K}x\inv(1-k)(1-h)x\overline{H} =
    \overline{K}(1-x\inv kx-x\inv hx+x\inv khx)\overline{H},$$
and using Lemma~\ref{trace} one has
    $$\matriz{{c}
    T(\overline{K}x\inv hx \overline{H})=T(\overline{K}x\inv kx \overline{H})=0, \\
    T(\overline{K}\,\overline{H})=|H|\quad \mbox{and} \quad T(\overline{K}x\inv kh x \overline{H})\ge 0.}$$
Therefore $T(ab)\ge |H| \ge 2=2T(1)$.  Then $(\beta_{x,h,H},\gamma_{x\inv,k,K})$ is a free pair by
Corollary~\ref{Salwa}. Let now $a=\beta_{x,h,H}-1$ and $b=\beta_{x\inv,xkx\inv,xKx\inv}-1$. Then
    $$\matriz{{rcl}
    ab&=&(1-h)x\overline{H}(1-xkx\inv)x\inv\overline{xKx\inv} \\ &=&
    x\overline{H}x\inv\overline{xKx\inv} - hx\overline{H}x\inv\overline{xKx\inv} - x\overline{H}xkx^{-2}\overline{xKx\inv} +
    hx\overline{H}xkx^{-2}\overline{xKx\inv}.}$$
As in the previous case $T(x\overline{H}x\inv\overline{xKx\inv})=|H|$ and
$T(hx\overline{H}xkx^{-2}\overline{xKx\inv})\ge 0$. Furthermore $(hx)\inv \not\in Hx\inv(xKx\inv)=HKx\inv=Kx\inv$
because $x\inv h x \not\in K$ and $x\inv \not\in Hxkx^{-2} (xKx\inv) = Hxkx\inv Kx\inv$, because $xkx\inv\not\in K=HK$.
Thus $1\not\in Hxkx\inv K$, hence Lemma~\ref{trace} yields $0=T(hx\overline{H}x\inv\overline{xKx\inv}) =
T(x\overline{H}xkx^{-2}\overline{xKx\inv})$ and so $T(ab)\ge 2$. Again one deduces that
$(\beta_{x,h,H},\beta_{x\inv,xkx\inv,xKx\inv})$ is a free pair from Corollary~\ref{Salwa}.
\end{proof}

Theorem~\ref{FreeComp} is a direct consequence of the following.

\begin{corollary}\label{FreeCompExp}
Let $x$ and $h$ be elements of a group and assume that $h$ has finite order. Then the following conditions are
equivalent.

\begin{tabular}{lll}
{\rm (1)} $x\inv h x \not\in \GEN{h}$. &\hspace{2cm} & {\rm (5)} $(\beta_{x,h},\beta_{x\inv,xhx\inv})$ is a free pair. \\
{\rm (2)} $x h x\inv \not\in \GEN{h}$. && {\rm (6)} $(\beta_{x,h},\gamma_{x\inv,h})$ is a free pair.  \\
{\rm (3)} $\beta_{x,h}\ne 1$. && {\rm (7)} $(\gamma_{x,h},\beta_{x\inv,h})$ is a free pair. \\
{\rm (4)} $\gamma_{x,h}\ne 1$. && {\rm (8)} $(\gamma_{x,h},\gamma_{x\inv,xhx\inv})$ is a free pair.
\end{tabular}

\end{corollary}

\begin{proof}
The equivalence between the first four conditions is obvious. The equivalence with (5) and (6) is a consequence of
Proposition~\ref{ManyFP} and the equivalence with (7) and (8) follows by symmetry.
\end{proof}

\section{Examples and questions. Bicyclic units}\label{ExBic}

Corollary~\ref{Nil} raised some natural questions. For a finite group $G$ let
    $$\matriz{{llll}
    M(G) & = & \min\{m\in \N : & (\alpha,\beta^m) \mbox{ is a free pair for every }
        \alpha=1+a \mbox{ and } \beta=1+b \\ & & & \mbox{ bicyclic units of } \Z G \mbox{ with } ab \mbox{ not nilpotent}\}, \\
    m(G) & = & \min\{m\in \N : & (\alpha,\beta^m) \mbox{ is a free pair for every }
        \alpha=1+a \mbox{ and } \beta=1+b \\ & & & \mbox{ bicyclic units of the same type of } \Z G
        \mbox{ with } ab \mbox{ not nilpotent}\}.}$$
By Theorem~\ref{Nil} one has: $G$ is Hamiltonian if and only if $M(G)=\infty$ if and only if $m(G)=\infty$. (Here we
are using the convention that the minimum of the empty set is $\infty$.)

\begin{problem}
Compute $M(G)$ and $m(G)$ for some (non-Hamiltonian) finite groups $G$.
\end{problem}

Finding a finite group $G$ with $1<M(G)\ne \infty$ is easy.

\begin{example}\label{S3}
$M(S_3)=2$.
\end{example}

\begin{proof}
Let $S_3$ be the symmetric group on three symbols. Recall that $S_3$ has two linear characters $\chi_1$ and $\chi_2$
and one irreducible character $\chi_3$ of degree $2$. Then $6T=\chi_1+\chi_2+2\chi_3$ is the character afforded by the
regular representation of $G$. If $\alpha=1+a$ and $\beta=1+b$ are bicyclic units then $\chi_1(ab)=\chi_2(ab)=0$ and
therefore $\chi_3(ab)=3 T(ab)$. Then applying Corollary~\ref{Char} one has: If $T(ab)=0$ then $\chi_3(ab)=0$ and
therefore $ab$ is nilpotent and $\GEN{\alpha,\beta}$ is a nilpotent group; if $|T(ab)|>1$ then $|\chi_3(ab)|\ge 4$ and
therefore $(\alpha,\beta)$ is a free pair. Otherwise, (that is if $T(ab)=\pm 1$), then $(\alpha,\beta)$ is not a free
pair (because $3$ is non-free) but $(\alpha,\beta^2)$ is a free pair. This shows that $M(S_3)\le 2$. To show that the
equality holds we exhibit a pair of bicyclic units $\alpha$ and $\beta$ such that $T(ab)=1$.

Let $\sigma=(1,2,3)$, $\tau=(1,2)$ and consider the bicyclic units $\alpha=\gamma_{\sigma,\tau}=1+a$ and
$\beta=\beta_{\sigma^2,\sigma\tau}=1+b$. We see that $T(ab)=1$ by using Lemma~\ref{trace} and noticing that
$1,\tau,\sigma^2\in \GEN{\tau}\GEN{\sigma\tau}=\{1,\tau,\sigma\tau,\sigma^2\}$,
$\sigma^2\tau\not\in\GEN{\tau}\GEN{\sigma\tau}$, $\GEN{\tau}\cap \sigma^2\GEN{\sigma\tau}\sigma^{-2} = \GEN{\tau}\cap
\GEN{\sigma\tau}=\{1\}$, $\GEN{\tau}\cap \tau\GEN{\sigma\tau}\tau\inv = \GEN{\tau}\cap \GEN{\sigma^2\tau}=\{1\}$ and
    $$ab=\overline{\tau}(1-\sigma\tau\sigma^2-\sigma\sigma\tau\sigma^2+\sigma\tau\sigma\tau\sigma^2)\overline{\sigma\tau} =
    \overline{\tau}(1-\sigma^2\tau-\tau+\sigma^2)\overline{\sigma\tau}.$$
\end{proof}

Finding a finite group $G$ such that $1<m(G)\ne \infty$ is more difficult. For example, if $D_{2n}$ denotes the
dihedral group of order $2n$ and $n$ is prime then $m(D_{2n})=1$ \cite{JRR}. This has been recently generalized by
Jim\'{e}nez \cite{J} who, using Corollary~\ref{Char}, have shown that $m(D_{2n})\le 2$ for every $n$ and, if $12$ does not
divide $n$, then $m(D_{2n})=1$. In fact the result of Jim\'{e}nez shows that if $12$ divides $n$ then $m(D_{2n})=1$ if and
only if $2\sqrt{3}$ is free, and  otherwise $m(D_{2n})=2$. More examples of finite groups $G$ with $m(G)=1$ can be
founded in \cite{DJR}. Since we still do not know if $2\sqrt{3}$ is free it is not clear whether dihedral groups
provides examples of non Hamiltonian groups $G$ for which $m(G)>1$. We found that $S_4$ provides such an example as an
application of Corollary~\ref{Char}.

\begin{example}
$m(S_4)=2$.
\end{example}

\begin{proof}
Having in mind that each irreducible character of $S_4$ takes values on the integers and that an integer $m$ is free if
and only if $|m|\ge 4$, it is not difficult to compute $m(S_4)$ using Corollary~\ref{Char} and an algebraic software
package. First, compute all the bicyclic units of one type: there are 157. Then compute $\chi(ab)$ for all the pairs
$(\alpha=1+a,\beta=1+b)$ of bicyclic units. It turns out that either $\chi(ab)=0$ for all the irreducible characters
$\chi$ of $S_4$ or there is an irreducible character $\chi$ such that $|\chi(ab)|\ge 2$. In the former case
$\GEN{\alpha,\beta}$ is nilpotent and in the latter $(\alpha,\beta^2)$ is a free pair. This shows that $m(S_4)\le 2$
and in fact the equality holds because there is a pair $(\alpha=1+a,\beta=1+b)$ such that $|\chi(ab)|$ is either $0$ or
$2$ for each irreducible character $\chi$. For the readers convenience we give a computer-free proof.

We choose the following set of generators of $S_4$.
        $$\sigma=(1,2,3), \quad \tau=(1,2), \quad \mu=(1,2)(3,4), \quad \nu=(1,3)(2,4).$$
Let $\pi:\Z S_4 \rightarrow \Z S_3$ be the ring homomorphism which acts as the identity on $S_3$ and maps $\mu$ and
$\nu$ to $1$. Let $\chi_1$, $\chi_2$ and $\chi_3$ be the three irreducible characters of $S_3$ (Example~\ref{S3}). Then
$S_4$ has two linear characters, $\theta_1=\chi_1\circ\pi$ and $\theta_2=\chi_2\circ\pi$; one irreducible character of
degree $2$, $\theta_3=\chi_3\circ\pi$; and two irreducible characters $\theta_4$ and $\theta_5$ of degree $3$. Observe
also that $\theta_i(g)\in \Z$ for each $g\in S_4$ and $i\le 5$. Furthermore $\theta_1(g)\equiv \theta_2(g)\equiv 1 \mod
2$, $\theta_4(g)\equiv \theta_5(g)\not\equiv \theta_3(g) \mod 2$ and if $\theta_3(g)$ is odd then $g$ is a $3$-cycle.

Let $\alpha=1+a$ and $\beta=1+b$ be bicyclic units of the same type of $\Z S_4$ such that $ab$ is not nilpotent but
$(\alpha,\beta)$ is not a free pair. Then $\pi(\alpha)$ and $\pi(\beta)$ are powers of bicyclic units of $\Z S_3$ and
$(\pi(\alpha),\pi(\beta))$ is not a free pair. Since $m(S_3)=1$, $\pi(ab)$ is nilpotent and thus
$0=\theta_1(ab)=\theta_2(ab)=\theta_3(ab)$. Therefore the number of 3-cycles in the support of $ab$ is even and thus
$\theta_4(ab)\equiv \theta_5(ab) \equiv \theta_1(ab)= 0 \mod 2$. Since $ab$ is not nilpotent, either $\theta_4(ab)\ne
0$ or $\theta_5(ab)\ne 0$ and thus the absolute value of one of them is at least $2$. Then either $2\theta_4(ab)$ or
$2\theta_5(ab)$ is free and so $(\alpha,\beta^2)$ is a free pair. This proves that $m(S_4)\le 2$.

To show $m(S_4)=2$ we consider the bicyclic units
    $$\alpha=\beta_{\mu\tau,\nu}=1+a \quad \mbox{and} \quad \beta=\beta_{\mu\nu\sigma^2\tau,\nu\sigma}=1+b.$$
As in Example~\ref{S3} we prove that $(\alpha,\beta)$ is not a free pair by showing that $\theta(x)$ is non-free for
every irreducible character $\chi$ of $S_4$. Obviously $\theta_1(ab)=\theta_2(ab)=0$. Furthermore $\pi(b)=0$ because
$\pi(\nu\sigma)=\sigma$ generates a normal subgroup of $S_3$, and thus $\theta_3(ab)=0$. Now a straightforward
computation shows that $\theta_4(ab)=-\theta_5(ab)=\pm 2$, which is non-free.
\end{proof}

The previous computations raised the following.

\begin{question}
Is there a positive integer $m$ such that $m(G)\le m$, for every non-Hamiltonian finite group $G$? Is $m(G)\le 2$ for
every non-Hamiltonian finite group $G$?
\end{question}

\section{Proof of Theorem~\ref{BB}}

The proof of Theorem~\ref{BB} uses the following Theorem of Gon\c{c}alves and Passman \cite{GPJA}.

\begin{theorem}\label{STau}
Let $V$ be a finite dimensional vector space $V$ over $\C$ and $S$ and $\tau$ endomorphisms of $V$. Assume that
$\tau^2=0$ and $S$ is diagonalizable. Let $r_+$ and $r_-$ be the maximum and minimum of the absolute values of the
eigenvalues of $S$. Let $V_+$ (resp. $V_-$) be the subspace generated by the eigenvectors of $V$ with eigenvalue of
modulus $r_+$ (resp. $r_-$) and $V_0$ the subspace generated by the remaining eigenvectors.

If the four intersections $V_{\pm} \cap \ker(\tau)$ and $\Imagen(\tau) \cap (V_0\oplus V_{\pm})$ are trivial then
$(S^s,(1+\tau)^t)$ is the free product of $\GEN{S^s}$ and $\GEN{(1+\tau)^t}$ for sufficiently large positive integers
$s$ and $t$.
\end{theorem}

Let $d$, $k$ and $m$ be positive integers such that $k$ and $d$ are relatively prime, and $m$ a multiple of $\phi(d)$.
Here $\phi$ stands for the Euler function. Since $k^m\equiv 1 \mod d$ the polinomial
    $$u_{k,m,d}(X)=(1+X+X^2+\ldots+X^{k-1})^m + \frac{1-k^m}{d}(1+X+X^2+\ldots+X^{d-1})$$
has integral coefficients.

Notice that the Bass cyclic units of $\Z G$ are the elements of the form $u_{k,m}(x)=u_{k,m,d}(x)$, where $x$ is an
element of order $d$ in the group $G$ and $k$ and $m$ satisfy the above conditions. If $\xi$ is a complex root $d$-th
root of unity then $u_{k,m,d}(\xi)=u_{k,m}(\xi)$ is a well defined unit of $\Z[\xi]$. Using this it is easy to see that
the Bass cyclic units of $\Z G$ belong to $U(\Z G)$. See \cite{GPJA} and \cite{Sehgal} for other properties of Bass
cyclic units. Notice that $u_{k,m}(x)$ is determined by $k$ modulo $d$, and hence one can assume that $1\le k \le d-1$.
Moreover $u_{1,m}(x)=1$, $u_{d-1,m}(x)=x^{(d-1)m}$ and $u_{k,m}(x)^a = u_{k,am}(x)$, for each integer $a$ \cite[Lemma
3.1]{GPJA}. Therefore the set of Bass cyclic units of $\Z G$ is closed under taking powers and $u_{k,m}(x)$ has finite
order if and only if $k\equiv \pm 1 \mod d$.

Let $G$ be a finite group of order coprime with $6$. We have to show that $\Z G$ has a free pair formed by a bicyclic
unit and a Bass cyclic unit. We start with a reduction argument.

\begin{quote}
{\em Claim 1. One may assume that the group $G=A\rtimes X$ where $X=\GEN{x}$ has prime order (say $p$), and one of the
following conditions hold:
\begin{enumerate}
\item $A$ is cyclic of prime power order.
\item $A$ is an elementary abelian $p$-group of order $p^2$.
\item $A$ is an elementary abelian $q$-group, with $q$ prime distinct from $p$, and $X$ acts faithfully and irreducibly
on $A$.
\end{enumerate}}
\end{quote}

Notice that if Theorem~\ref{BB} holds for some subgroup or some epimorphic image of $G$ then it also holds for $G$.
This is clear if $H$ is a subgroup of $G$ because a Bass cyclic unit (respectively, bicyclic unit) of $\Z H$ is also a
Bass cyclic unit (respectively, bicyclic unit) of $\Z G$. Assume that $H$ is an epimorphic image of $G$, $u_1$ is a
Bass cyclic unit of $\Z H$ and $\beta_1$ is a bicyclic unit of $\Z H$ such that $(u_1,\beta_1^s)$ is free for $s$
sufficiently large. Then $\Z G$ has a Bass cyclic unit $u$ and a bicyclic unit $\beta$ such that $u$ projects to
$u_1^k$ and $\beta$ projects to $\beta_1^l$ for some $k,l\ge 0$ \cite[Lemma 3.2]{GPJA}. Thus $(u,\beta^s)$ is a free
pair of $\Z G$, for $s$ sufficiently large. Thus arguing by induction to prove Theorem~\ref{BB}, one may assume without
loss of generality that all the proper subgroups and proper epimorphic images of $G$ are abelian. Then a theorem of
\cite{MM} shows that the group $G$ is as stated in the claim. This proves Claim 1.

So, in the remainder we assume that $G$ is as in Claim 1, with $p,q\ge 5$, and we refer to cases (1), (2) and (3),
depending on which condition holds. Since $G$ is non-abelian, $X$ is not normal in $G$, and therefore $x^a\not\in X$
for some $a\in A$ that will be fixed until the end of the proof. For example, in Case (1), $a$ can be a generator of
$A$. In Case (2) one can take $a\in A$ such that $Z(G)=\GEN{b=a^xa\inv}$. In Case (3), $Z(G)=1$, and thus every
non-trivial element of $A$ satisfies the required condition.

Let $x$ and $a$ be as above and let $q$ be the order of $a$. (Notice that this notation is compatible with the notation
in Case (3); in Case (1) $q$ is a power of $p$ and; in Case (2) $q=p$.) We fix $2\le k \le q-2$ (recall that $q\ge 5$)
and let
    $$u=u_{k,m}(a), \quad \mbox{and} \quad \beta=1+(1-x)a\widehat{x}$$
for $m$ a multiple of $\phi(d)$. Later on we will specify some additional condition on $m$. Notice that $u$ and $\beta$
have infinite order.

\begin{quote}
{\em Claim 2. There is a linear representation $\chi$ of $A$ such that the induced representation $\rho=\chi^G$ is
irreducible, $\rho((a,x))\ne 1$ and in Case (3), either $|A|=q$ or $a\in \ker(\chi)$.}
\end{quote}

Since the commutator $(a,x)\ne 1$, there is an irreducible representation (necessarily non-linear) $\rho$ of $G$ such
that $\rho((a,x))\ne 1$. All the non-linear irreducible representations of $G$ are induced from linear representations
of $A$ and so the claim is clear except for the exceptional case (3) with $|A|\ne q$. In this case $A$ has a maximal
subgroup $B$ containing $a$ and so there is a linear representation $\chi$ of $A$ with $\ker(\chi)=B$. Since $X$ acts
irreducibly on $A$ and the subgroups $G'$ and $C=\GEN{a^{x^i} : i=0,1,\dots,p-1}$ of $A$ are invariant under this
action, one has $A=G'=C$. The first equality implies that $\chi^G$ is irreducible and the second that
$a^{x^i}\not\in\ker(\chi)$ for some $i$. Therefore
$\rho(a)=\diag(\chi(a)=1,\chi(a^x),\chi(a^{x^2}),\dots,\chi(a^{x^{p-1}})) \ne
\diag(\chi(a^x),\chi(a^{x^2}),\dots,\chi(a^{x^{p-1}}),\chi(a)) = \rho(a^x)$ and thus $\rho((a,x))\ne 1$, as wanted.
This proves Claim 2.

For the remainder of the proof we fix a linear representation $\chi$ of $A$ satisfying the conditions of Claim 2 and
put $\zeta_i=\chi(a^{x^i})$ and $\rho=\chi^G$. Then $\rho(a)=\diag(\zeta_0,\zeta_1,\dots,\zeta_{p-1})$ is a non-scalar
matrix, i.e. the set $\Lambda=\{\zeta_i:i=0,1,\dots,p-1\}$ has at least two different elements. In fact

\begin{quote}
{\em Claim 3. If $|\Lambda|\ne p$ then Case (3) holds and $|A|\ne q$ and so $\Lambda$ contains $1$ (and another
different element).

Furthermore, in Case (1) $\zeta_i$ and $\zeta_{i+1}$ are not complex conjugate for each $0\le i<p$ (where the
subindexes are considered as elements in $\Z_p$, the ring of integers modulo $p$).}
\end{quote}

We consider the three cases separately.

In Case (1), $a^x=a^r$ for some integer $r$ coprime with $p$ and such that the multiplicative order of $r$ modulo the
order of $a$ is $p$. Therefore $\zeta_i=\zeta_0^{r^i}$ and $\zeta_0\mapsto \zeta_0^r$ induces an automorphism $\sigma$
of $\Q(\zeta_0)$ ($=\Q(\zeta_i)$) of order $p$. This implies that $|\Lambda|=p$. Moreover, if $\zeta_i$ and
$\zeta_{i+1}$ are complex conjugate then $\sigma$ has order $2$ yielding a contradiction with $p \ge 5$.

In Case (2), $Z(G)=G'=\GEN{b=a^xa\inv}$. Then $\zeta_i=\xi^i\zeta_0$, where $\xi=\chi(b)$ is a primitive $p$-th root of
unity, and again $|\Lambda|=p$.

In Case (3) with $|A|\ne q$, $1=\chi(a)=\zeta_0 \in \Lambda$, by the construction of $\chi$. If $|A|=q$ then $\chi$ is
injective and thus the $\zeta_i$'s are pairwise different. Indeed, if $0\le i < j \le p-1$ then $x^{j-i}$ is a
generator of $X$. Then $a^{x^{j-i}}\ne a$, so $\zeta_i=\chi(a^{x^i})\ne \chi(a^{x^j})=\zeta_j$. This finishes the proof
of Claim 3.

The representations of $u$ and $\beta$ by $\rho$ are
    $$S=\rho(u)=\diag(u_{k,m}(\zeta_0),u_{k,m}(\zeta_1),\ldots,u_{k,m}(\zeta_{p-1}))=\diag(u_0,u_1,\dots,u_{p-1})$$
and $\rho(\beta)=1+\tau$ with
    $$\tau = \pmatriz{{cccc}
    \zeta_0-\zeta_1 & \zeta_0-\zeta_1 & \dots & \zeta_0-\zeta_1 \\
    \zeta_1-\zeta_2 & \zeta_1-\zeta_2 & \dots & \zeta_1-\zeta_2 \\
    \vdots & \vdots &  & \vdots \\
    \zeta_{p-1}-\zeta_0 & \zeta_{p-1}-\zeta_0 & \dots & \zeta_{p-1}-\zeta_0}.
    $$
Since not all the $\zeta_i$ are equal, $\tau$ has rank $1$ with image generated by
    $$\Psi=\pmatriz{{c} \zeta_0-\zeta_1 \\ \zeta_1-\zeta_2 \\\vdots \\
    \zeta_{p-1}-\zeta_0},$$
and kernel
    $$K=\ker(\tau) = \left\{\pmatriz{{c}x_0\\\vdots\\x_{p-1}}: x_0+\cdots+x_{p-1}=0\right\}.$$

Let $r_+$ and $r_-$ be the maximum and minimum of $\{|u_i|:i=0,1,\dots,p-1\}$ and set $X_+=\{i:|u_i|=r_+\}$,
$X_-=\{i:|u_i|=r_-\}$ and $X_0=\Z_p\setminus (X_+\cup X_-)$. Let $\{e_1,\dots,e_n\}$ be the canonical basis of
$V=\C^p$, the representation space of $\rho$. Let $V_*$ be the span $\{e_i:i\in X_*\}$ for $*=+,-$ or $0$. Notice that
this notation agrees with that of Theorem~\ref{STau}.

\begin{quote}
{\em Claim 4: $\Z_p\ne X_+$ (equivalently $\Z_p\ne X_-$).}
\end{quote}

The order $d$ of $\rho(a)$ is a divisor of $q$ and so it is a prime power. Moreover $\Lambda\subseteq L$ where $L$ is
the set of $d$-th roots of unity. By \cite[Lemma 3.5(ii)]{GPJA}, if $x,y\in L$ then $|u_{k,m}(x)|=|u_{k,m}(y)|$ if and
only if $x$ and $y$ are either equal or conjugate. If $\Z_p=X_+$ then $\Lambda$ has either one element or it has two
different conjugate elements. This contradicts Claim 3, and so $\Z_p\ne X_+$. This proves Claim 4.

At this point it is tempting to try to show that $S$ and $\tau$ satisfy the conditions of Theorem~\ref{STau}.
Unfortunately this is not the case in general. For example, in Case (3) some of the $\zeta_i$'s may be equal (for
example, this is the case if $q<p$) and this may provide some non-trivial elements in $V_+\cap K$.

As we mentioned above, we are going to be more specific on the the integer $m$ used in the definition of
$u=u_{k,m}(a)$. Namely, we are going to impose that $m$ is a multiple of $q$, the order of $a$. Let $\xi$ be a
primitive $q$-th root of unity. If $|u_{k,m}(\xi^a)|=|u_{k,m}(\xi^b)|$ then $b\equiv \pm a \mod q$ (\cite[Lemma
3.5]{GPJA}). If moreover $a\not\equiv 0 \mod q$ then $u_{k,m}(\xi^b)=\left(\frac{\xi^{bk}-1}{\xi^b-1}\right)^m
=\left(\frac{\xi^{bk}}{\xi^b} \cdot \frac{\xi^{ak}-1}{\xi^a-1}\right)^m =\left(\frac{\xi^{ak}-1}{\xi^a-1}\right)^m =
u_{k,m}(\xi^a)$. In particular, $\{\zeta_i:i\in X_+\}$ and $\{\zeta_i:i\in X_-\}$ have cardinality 1. This implies that
$W=(V_+\cap K) \oplus (V_-\cap K)$ is invariant under the action of $S$ and hence the endomorphisms $S$ and $\tau$ of
$V$ induce endomorphisms $\overline{S}$ and $\overline{\tau}$ of $\overline{V}=V/W$.

We are going to use the standard bar notation for reduction modulo $W$. Then, as we will see below,
$\overline{V}=\overline{V_+}\oplus \overline{V_0} \oplus \overline{V_-}$ is a decomposition of $\overline{V}$ with
respect to $\overline{S}$ as in Theorem~\ref{STau}.

\begin{quote}
{\em Claim 5: $\overline{S}$ and $\overline{\tau}$ satisfy the hypothesis of Theorem~\ref{STau}.}
\end{quote}

The kernel and image of $\overline{\tau}$ are $K_1=\overline{\tau\inv(W)}$ and $I_1=\overline{I}$, respectively. The
dimension of $\overline{V_{\pm}}\simeq V_{\pm}/(V_{\pm} \cap W)= V_{\pm}/(V_{\pm} \cap K)$ is 1, because $K$ is a hyperplane of $V$ but
$e_i\not\in K$ and $S(e_i)=\zeta_i e_i \not\in K$ for each $i\in X_{\pm}$. This shows that $\overline{V_{\pm}} \cap K_1 = 0$.

To prove that $I_1 \cap (\overline{V_0} \oplus \overline{V_{\pm}})=0$ we only have to show that $\overline{\Psi}\not\in
\overline{V_0} \oplus \overline{V_{\pm}}$, and due to symmetry, only that $\overline{\Psi}\not\in\overline{V_0} \oplus
\overline{V_+}$. If this is not so, then $\Psi = v + w_+ +w_-$ for some $v\in V_0$, $w_+\in V_+$ and $w_-\in V_-\cap
K$. Thus the projection of $\Psi$ on $V_-$ belongs to $K$. That is, $\sum_{i\in X_-}(\lambda-\zeta_{i+1}) = \sum_{i\in
X_-} (\zeta_i-\zeta_{i+1})=0$, where $\lambda$ is the unique element of the form $\zeta_i$, with $i\in X_-$. Now we
delete from the equality above the neighbor zero summands, that is the ones in which $\lambda=\zeta_{i+1}$. Passing the
negative expressions to the right side and adding up the equal terms, we obtain
    $n\lambda = \sum_{i=1}^k m_i \mu_i$,
for $n$ the cardinalty of $\{i\in X_-:\zeta_{i+1}\ne \lambda\}$, $m_i$ non-negative integers and $\mu_1,\dots,\mu_i$,
$q$-th roots of unity different from $\lambda$. Moreover $n>0$, because $\Z_p\ne X_-$. Let $\tr$ denotes the Galois
trace of the extension $\Q(\xi)/\Q$ and write $q=p^r$, with $p$ prime. Let $\varepsilon$ be a $q$-th root of unity.
Then $\tr(\varepsilon)=p^{r-1}(p-1)$ if $\varepsilon=1$, $\tr(\varepsilon)=-p^{r-1}$ if $\varepsilon$ has order $p$ and
$\tr(\varepsilon)=0$ otherwise. Thus $\tr(\lambda\inv \mu_i)\le 0$ for each $i$ and hence
$0<n(p-1)p^{n-1}=\tr(n)=\sum_{i=1}^k m_i \tr(\lambda\inv \mu_i)\le 0$, which yields a contradiction. This finishes the
proof of Claim 5.\bigskip

By Claim 5, $(\overline{S}^s,(1+\overline{\tau})^t)$ is the free product of $\GEN{\overline{S}^s}$ and
$\GEN{(1+\overline{\tau})^t}$ for $s$ and $t$ sufficiently large. Thus $(u^s,\beta^t)$ the free product of $\GEN{u^s}$
and $\GEN{\beta^t}$. Since $u$ and $\beta$ have infinite order, $(u^s,\beta^t)$ is a free pair formed by a Bass cyclic
unit and a power of a bicyclic unit of $\Z G$. This finishes the proof of Theorem~\ref{BB}.

\section{Examples and questions. Bass cyclic and bicyclic units}

If $u_{k,m}(a)$ is a Bass cyclic unit of infinite order then $k\not \equiv \pm 1 \mod d$, where $d$ is the order of $a$
and hence $d$ is not a divisor of $4$ nor $6$. Therefore, if $G$ has a free pair in which one of the elements is a Bass
cyclic unit, and the other is a power of a bicyclic unit, then $G$ has a non central element whose order does not
divide neither $4$ nor $6$ (for $\Z G$ to have a non-central Bass cyclic unit of infinite order) and $G$ is not
Hamiltonian (for $\Z G$ to have a non-trivial bicyclic unit). This justifies partially the hypothesis in
Theorem~\ref{BB} (and in \cite[Theorem 4.7]{GPJA}) of $G$ having order coprime to 6 and suggest the following question.

\begin{question}\label{QBB}
Let $G$ be a non-Hamiltonian group with an element whose order does not divide neither $4$ nor $6$. Does $\Z G$ have a
free pair formed by a Bass cyclic and a power of a bicyclic unit?
\end{question}

\begin{remark}
Notice that if $G$ is non-abelian and has a central element $a$ of order $d$ with $d\nmid 4$ and $d\nmid 6$, then $G$
has also a non-central element satisfying the same condition. Indeed, assume that every non-central element of $G$ has
order dividing either $4$ or $6$. Let $b$ a non-central element of maximal order $n$, and so $n=2,3,4$ or $6$. Then
$a^ib$ is non-central and therefore the order of $a^ib$ divides $4$ or $6$ for each $i$. This implies that $\GEN{a}\cap
\GEN{b}\ne 1$ and thus $n\ne 2,3$. Furthermore $d=2^k 3^l$ for some $k,l\ge 0$. If $n=6$ then the orders of $a^{3^l}b$
and $a^{2^k}b$ are multiple of $2^k \cdot 3$ and $2\cdot 3^l$, respectively, and this yields a contradiction because
either $k\ge 2$ or $l\ge 2$. This shows that the order of $b$ is $4$ and the order of $ab$ is at most $4$. Then
$b^3\not\in \GEN{a}$ and so the order of $ab$ divides $4$. Thus the order of $a$ divides $4$, contradicting the
hypothesis.
\end{remark}

Notice that if either $A_n$ or $S_n$ satisfies the hypothesis of Question~\ref{QBB} then $n\ge 5$. Thus to give an
affirmative answer to Question~\ref{QBB} for symmetric and alternating groups, it is enough to show that $\Z A_5$ has a
free pair formed by a power of a  bicyclic unit, and a Bass cyclic unit.

The group $A_5$ can be defined by generators and relations as
    $$A_5=\GEN{a,b|a^2=b^3=(ba)^5=1}$$
and one can take $a=(1,2)(3,4)$ and $b=(1,3,5)$, so that $c=ba=(1,2,3,4,5)$. Take the irreducible representation $\phi$
of $A_5$ given by
    $$A=\phi(a)=\pmatriz{{rrrr} 0 & 1 & 0 & 0 \\ 1 & 0 & 0 & 0 \\ 0 & 0 & 0 & 1 \\ 0 & 0 & 1 & 0} \quad \mbox{and}
    \quad
    B=\phi(b)=\pmatriz{{rrrr} 0 & 0 & 1 & 0 \\ 0 & 1 & 0 & 0 \\ -1 & -1 & -1 & -1 \\ 0 & 0 & 0 & 1}.$$
Let $\xi=e^{\frac{2\pi i}{5}}$, a primitive 5-th root of unity, and set $F=\Q(\xi)$. Let $\sigma\in \Gal(F/\Q)$ be
given by ($\sigma(\xi)=\xi^2$. Notice that $\sigma$ generates $\Gal(F/\Q)$ and one has
    $$C=\phi(c)=BA = \pmatriz{{rrrr} 0 & 0 & 0 & 1\\1 & 0 & 0 & 0\\-1 & -1 & -1 & -1\\0 & 0 & 1 & 0}.$$
Consider $\Gal(F/\Q)$ acting componentwise on $F^4$. Let $v_0=(\xi^2,\xi,\xi^4,\xi^3)$ and $v_i=\sigma^i(v_1)$. Then
$Cv_0=\xi v_0$, that is $v_0$ is an eigenvector of $C$ with eigenvalue $\xi$. Therefore $v_i$ is also an eigenvector of
$C$ with eigenvalue $\sigma^i(\xi)=\xi^{2^i}$. This implies that $v_i$ is an eigenvector of the Bass cyclic unit
$S=u_{2,4}(C)=\phi(u_{2,4}(c))$ with eigenvalue $\lambda_i=\sigma^i(u_{2,4}(\xi))=1+\xi^{2^i}$. Now
$|\lambda_1|=|\lambda_3|>|\lambda_2|=|\lambda_4|$ and so $F^4 = V_+\oplus V_-$, where $V_+=\GEN{v_1,v_3}$ and
$V_-=\GEN{v_2,v_4}$.

Let now
    $$\tau = \phi((1-a)b(1+a)) =
    \pmatriz{{rrrr} -1 & -1 & 1 & 1\\1 & 1 & -1 & -1\\-2 & -2 & -3 & -3\\2 & 2 & 3 & 3}.$$
Then $K=\ker(\tau) = \Imagen(\tau)= \{(x_1,x_2,x_3,x_4) : x_1+x_2=x_3+x_4=0\}$. Let $\delta_1,\delta_2\in F$ be such
that $\delta_1 v_1 + \delta_2 v_3 = (\delta_1 \xi^2 + \delta_2 \xi^3,\delta_1 \xi + \delta_2 \xi^4 ,\delta_1 \xi^4 +
\delta_2 \xi ,\delta_1 \xi^3 + \delta_2 \xi^2) \in K$. Then $\delta_1 (\xi^2+\xi) + \delta_2 (\xi^3+\xi^4) =
\delta_1(\xi^4+\xi^3) + \delta_2(\xi+\xi^2) = 0$ and so $\delta_1=\delta_2=0$, because $(\xi^4+\xi^3)-(\xi^2+\xi) =
-2\xi^2-2\xi-1 \ne 0$. This shows that $K\cap X_+=0$. Since $\sigma$ leaves $K$ invariant and interchanges $V_+$ and
$V_-$ one obtains that $K\cap X_-=0$. Therefore $S$ and $\tau$ (considered as endomorphisms of complex vector spaces)
satisfy the hypotheses of Theorem~\ref{STau}. Thus $(S^n=\phi(u_{2,4}(c)^n),\phi(\beta_{a,b})^m =(1+\tau))^m$ is a free
pair for some $n$ and $m$, and we conclude that $(u_{2,4}(c)^n=u_{2,4n}(c),(\beta_{a,b})^m)$ is a free pair of $\Z
A_5$, formed by a Bass cyclic unit and a power of a bicyclic unit. So we have shown.

\begin{theorem}
$\Z A_n$ (resp. $\Z S_n$) contains a free pair formed by a power of a  bicyclic unit and a Bass cyclic unit, if and
only if $n\ge 5$.
\end{theorem}

{\bf Acknowledgements}. In this paper we have used many ideas which Donald Passman generously shared with us. We would
like to thank him and also thank the referee for pointing out an error in a previous version of the proof of
Theorem~\ref{BB}.


\end{document}